\title{Numerical semigroups, cyclotomic polynomials and Bernoulli numbers}
\author{Pieter Moree}
\documentclass[12pt]{article}
\usepackage{amssymb, latexsym, amsfonts}
\usepackage[T1]{fontenc}
\def\@ptsize{2}
\newtheorem{Thm}{Theorem}
\newtheorem{Problem}{Problem}

\newtheorem{Lem}{Lemma}

\newtheorem{cor}{Corollary}

\newcommand{\qed}{\hfill $\Box$}

\begin{document}
\date{}
\maketitle
{\def\thefootnote{}
}
\begin{abstract}
\noindent We give two proofs of a folklore result relating numerical semigroups 
of embedding dimension two and binary cyclotomic polynomials and explore some consequences. 
In particular, we give a more conceptual reproof of a result of Hong et al. (2012) on
gaps between the exponents of non-zero monomials in a binary cyclotomic polynomial.\\
\indent The intent of the author with this paper is to better unify  the various
results within the  cyclotomic polynomial and numerical semigroup communities. 
\end{abstract}
\section{Introduction}
Let $a_1,\ldots,a_m$ be positive integers, and let $S=S(a_1,\ldots,a_m)$ be the
set of all non-negative integer linear combinations of $a_1,\ldots,a_m$, that is,
$$S=\{x_1a_1+\cdots +x_m a_m~|~x_i\in \mathbb Z_{\ge 0}\}.$$
Then
$S$ is a {\it semigroup} (that is, it is closed under addition). The semigroup $S$ is
said to be {\it numerical} if its complement $\mathbb Z_{\ge 0}\backslash S$ is finite. 
It is not difficult to prove that $S(a_1,\ldots,a_m)$ is numerical if and only if $a_1,\ldots,a_m$ 
are relatively prime (see, e.g., \cite[p. 2]{RA}).
If $S$ is numerical,
then $\max \{\mathbb Z_{\ge 0}\backslash S\}=F(S)$ is the {\it Frobenius number} of $S$. 
Alternatively, by setting $d(k,a_1,\ldots,a_m)$ equal to the number of non-negative 
integer representations of $k$ by $a_1,\ldots,a_m$, one can characterize $F(S)$ as the largest
$k$ such that $d(k,a_1,\ldots,a_m)=0$. The value $d(k,a_1,\ldots,a_m)$ is called the 
{\it denumerant} of $k$. That $F(S(4,6,9,20))=11$ is well-known 
to fans of
Chicken McNuggets, as 11 is the largest number of McNuggets that cannot be 
exactly purchased; hence the notion of 
of the Frobenius number 
is less abstract than it
might appear at first glance.
A set of generators of a numerical semigroup is a minimal system of generators if none of its proper subsets generates the numerical semigroup. It is known
that every numerical semigroup $S$ has a unique minimal system of generators and also that this minimal system of generators is 
finite (see, e.g., \cite[Theorem 2.7]{roos}). The cardinality of the minimal set of generators is called the {\it embedding dimension} of the numerical semigroup $S$ and is denoted by $e(S)$. The smallest member in the minimal system of generators is called the {\it multiplicity} of the numerical semigroup $S$ and is denoted by $m(S)$.
The {\it Hilbert series} of the numerical semigroup $S$ is the formal
power series $$H_S(x)=\sum_{s\in S}x^s\in \mathbb Z[[x]].$$ 
It is practical to multiply this by $1-x$ as we then obtain a {\it polynomial}, called
the {\it semigroup polynomial}:
\begin{equation}
\label{semigppoly}
P_S(x)=(1-x)H_S(x)=x^{F(S)+1}+(1-x)\sum_{0\le s\le F(S)\atop s\in S}x^s=1+(x-1)\sum_{s\not\in S}x^s.
\end{equation}
{}From  $P_S$ one immediately reads off the
Frobenius number:
\begin{equation}
\label{frob}
{\rm deg}(P_S(x))=F(S)+1.
\end{equation}
\indent The $n$th cyclotomic polynomial $\Phi_n(x)$ is defined by
$$\Phi_n(x)=\prod_{1\le j\le n\atop (j,n)=1}(x-\zeta_n^j)=\sum_{k=0}^{\varphi(n)}a_n(k)x^k,$$ with
$\zeta_n$ a $n$th primitive root of unity (one can take $\zeta_n=e^{2\pi i/n}$).
It has degree $\varphi(n)$, with $\varphi$ Euler's totient function. The polynomial $\Phi_n(x)$ is irreducible over the rationals, see, e.g., Weintraub \cite{wein}, and has integer
coefficients. The polynomial $x^n-1$
factors as
\begin{equation}
\label{prod}
x^n-1=\prod_{d|n}\Phi_d(x)
\end{equation}
over the rationals. By M\"obius inversion it follows from (\ref{prod}) that
\begin{equation}
\label{prod1}
\Phi_n(x)=\prod_{d|n}(x^d-1)^{\mu(n/d)},
\end{equation}
where $\mu(n)$ denotes the M\"obius function. 
{}From (\ref{prod1}) one deduces that if $p|n$ is a prime, then
\begin{equation}
\label{pn}
\Phi_{pn}(x)=\Phi_n(x^p).
\end{equation}
A good source for further properties of cyclotomic polynomials is
Thangadurai \cite{Thanga}.\\
\indent A purpose of this paper is to popularise the following folklore result and point out
some of its consequences.
\begin{Thm}
\label{basic}
Let $p,q>1$ be coprime integers, then
$$P_{S(p,q)}(x)=(1-x)\sum_{s\in S(p,q)}x^s={(x^{pq}-1)(x-1)\over (x^p-1)(x^q-1)}.$$ 
\end{Thm}
In case $p$ and $q$ are distinct primes it follows from (\ref{prod1}) and Theorem \ref{basic} that
\begin{equation}
\label{corr}
P_{S(p,q)}(x)=\Phi_{pq}(x). 
\end{equation}
Already Carlitz \cite{Carlitz} in 1966 implicitly mentioned this result without proof.\\
\indent The Bernoulli numbers $B_n$ can be defined by
\begin{equation}
\label{generating}
{z\over e^z-1}=\sum_{n=0}^{\infty}B_n{z^n\over n!},~|z|<2\pi.
\end{equation}
One easily sees that $B_0=1,B_1=-1/2,B_2=1/6,B_3=0,B_4=-1/30$ 
and $B_n=0$ for all odd $n\ge 3$. 
The most basic recurrence relation is, for $n\ge 1$,
\begin{equation}
\label{basicrecurrence}
\sum_{j=0}^n {n+1\choose j}B_j=0.
\end{equation}
The Bernoulli numbers first arose in the
study of power sums $S_j(n):=\sum_{k=0}^{n-1}k^j$. Indeed, one has, cf. Rademacher \cite{Rademacher},
\begin{equation}
\label{powersum}
S_j(n)={1\over j+1}\sum_{i=0}^j {j+1\choose i}B_i n^{j+1-i}.
\end{equation}
In Section \ref{vijf}, we consider an infinite family of recurrences for $B_m$ of
which the following is typical
\begin{eqnarray*}
B_m & = & {m\over 4^m-1}(1+2^{m-1}+3^{m-1}+5^{m-1}+6^{m-1}+9^{m-1}+10^{m-1}+13^{m-1}+17^{m-1})\cr
& & +{7^m\over 4(1-4^m)}\sum_{r=0}^{m-1}{m\choose r}\Big({4\over 7}\Big)^r(1+2^{m-r}+3^{m-r})B_r.
\end{eqnarray*}
The natural numbers $1,2,3,5,6,9,10,13$ and $17$ are precisely those that are not in the 
numerical semigroup $S(4,7)$.\\
\indent Let $f=c_1x^{e_1}+\cdots+c_sx^{e_s},$ where the coefficients $c_i$ are non-zero and $e_1<e_2<\cdots<e_s$.
Then the {\it maximum gap} of $f$, written as $g(f)$, is defined by
$$g(f)=\max_{1\le i<s}(e_{i+1}-e_i),~g(f)=0{\rm ~when~}s=1.$$
Hong et al. \cite{zuidkorea} studied $g(\Phi_n)$ (inspired by a cryptographic application \cite{zuidkorea2}).
They reduce the study of these gaps to the case where $n$ is square-free and odd and established the following 
result for the simplest non-trivial case.
\begin{Thm} {\rm \cite{zuidkorea}.}
\label{hoon}
If $p$ and $q$ are arbitrary primes with $2<p<q$, then $g(\Phi_{pq})=p-1$. 
\end{Thm}
In Section \ref{zes} a conceptual proof of Theorem \ref{hoon} using numerical semigroups is given.
\section{Inclusion-exclusion polynomials}
It will turn out to be convenient to work with a generalisation of the cyclotomic polynomials, introduced
by Bachman \cite{B}.
Let $\rho=\{r_1,r_2,\ldots,r_s\}$ be a set of natural numbers satisfying $r_i>1$ and
$(r_i,r_j)=1$ for $i\ne j$, and put
$$n_0=\prod_i r_i,~n_i={n_0\over r_i},~n_{ij}={n_0\over r_ir_j}[i\ne j],\ldots~.$$
For each such $\rho$ we define a function $Q_{\rho}$ by
\begin{equation}
\label{ijs}
Q_{\rho}(x)={(x^{n_0}-1)\cdot \prod_{i<j}(x^{n_{ij}}-1)\cdots\over \prod_i (x^{n_i}-1)\cdot \prod_{i<j<k}(x^{n_{ijk}}-1)\cdots}.
\end{equation}
For example, if $\rho=\{p,q\}$, then
\begin{equation}
\label{prod4}
Q_{\{p,q\}}(x)={(x^{pq}-1)(x-1)\over (x^p-1)(x^q-1)}.
\end{equation}
It can be shown that $Q_{\rho}(x)$ defines a polynomial of degree $d:=\prod_i (r_i-1)$.
We define its coefficients $a_{\rho}(k)$ by
$Q_{\rho}(x)=\sum_{k\ge 0}a_{\rho}(k)x^k$. Furthermore, $Q_{\rho}(x)$ is {\it selfreciprocal}; that
is $a_{\rho}(k)=a_{\rho}(d-k)$ or, what amounts to the same thing, 
\begin{equation}
\label{self}
Q_{\rho}(x)=x^dQ_{\rho}({1\over x}).
\end{equation}
\indent If all elements of $\rho$ are prime, then comparison of (\ref{ijs}) with (\ref{prod1}) shows that
\begin{equation}
\label{qn}
Q_{\rho}(x)=\Phi_{r_1r_2\cdots r_s}(x).
\end{equation}
If $n$ is an arbitrary integer and $\gamma(n)=p_1\cdots p_s$ its squarefree kernel, then by (\ref{pn}) 
and (\ref{qn}) we have
$Q_{\{p_1,\ldots,p_s\}}(x^{n/\gamma(n)})=\Phi_n(x)$
and hence inclusion-exclusion polynomials generalize cyclotomic polynomials. They can be
expressed as products of cyclotomic polynomials.
\begin{Thm} \label{bach} {\rm \cite{B}}.
Given $\rho=\{r_1,\ldots,r_s\}$ and
$$D_{\rho}=\{d:d|\prod_i r_i{\rm ~and~}(d,r_i)>1{\rm ~for~all~}i\},$$
then $Q_{\rho}(x)=\prod_{d\in D_{\rho}}\Phi_d(x)$.
\end{Thm}
{\tt Example}. We have $Q_{\{4,7\}}=\Phi_{28}\Phi_{14}$.
\subsection{Binary inclusion-exclusion polynomials: a close-up}
Lam and Leung \cite{lam} discuss binary cyclotomic polynomials $\Phi_{pq}$
in detail, with $p$ and $q$ primes (their results were
anticipated by Lenstra \cite{Lenstra}). Now, let $p,q>1$ be positive coprime integers.
All arguments in their paper easily generalize to this setting (instead of taking
$\xi$ to be a primitive $pq$th-root of unity as they do, one has to take $\zeta$ a $pq$th root
of unity satisfying $\zeta^p\ne 1$ and $\zeta^q\ne 1$). One finds that
\begin{equation}
\label{leungie}
Q_{\{p,q\}}(x)=\sum_{i=0}^{\rho-1}x^{ip}\sum_{j=0}^{\sigma-1}x^{jq}-x^{-pq}
\sum_{i=\rho}^{q-1}x^{ip}\sum_{j=\sigma}^{p-1}x^{jq},
\end{equation}
where $\rho$ and $\sigma$ are the (unique) non-negative
integers for which $1+pq=\rho p+\sigma q$.
On noting that upon expanding the products in identity (\ref{leungie}), the resulting monomials
are all different, we arrive at the following result.
\begin{Lem}
\label{binary} 
Let $p,q>1$ be coprime integers. Let $\rho$ and $\sigma$ be the (unique) non-negative
integers for which $1+pq=\rho p+\sigma q$.
Let $0\le m<pq$. Then either $m=\alpha p+\beta q$ or $m=\alpha p+\beta q-pq$
with $0\le \alpha \le q-1$ the unique integer such that $\alpha  p\equiv m({\rm mod~}q)$
and $0\le \beta \le p-1$ the unique integer such that $\beta q\equiv m({\rm mod~}p)$.
The inclusion-exclusion coefficient $a_{\{p,q\}}(m)$ equals
$$\cases{1 & if $m=\alpha p+\beta q$ with $0\le \alpha\le \rho -1,~0\le \beta\le
\sigma -1$;\cr -1 & if $m=\alpha p+\beta q-pq$ with $\rho\le \alpha \le q-1,~\sigma\le 
\beta \le p-1$;\cr  0 & otherwise.}
$$
\end{Lem}
\begin{cor}
\label{carla}
The number of positive coefficients in $Q_{\{p,q\}}(x)$ equals $\rho \sigma$ and the number
of negative ones equals $\rho \sigma -1$. The number of non-zero coefficients equals $2\rho\sigma-1$.
\end{cor}
This corollary (in case $p$ and $q$ are distinct primes) is due
to Carlitz \cite{Carlitz}.

Lemma \ref{binary} can be nicely illustrated with an LLL-diagram (for Lenstra, Lam and Leung). 
Here is one such diagram for $p=5$ and $q=7$.
\begin{center}\begin{tabular}{ccc|cccc}
28&33&3&8&13&18&23\\
21&26&31&1&6&11&16\\\hline
14&19&24&29&34&4&9\\
7&12&17&22&27&32&2\\
0&5&10&15&20&25&30
\end{tabular}\end{center}
We start with 0 in the lower left and add $p$ for every move to the right and $q$ for every move upwards. Reduce modulo $pq$. Every integer $0,\ldots,pq-1$ is obtained precisely once
in this way (by the Chinese remainder theorem).\\ 
\indent Lemma \ref{binary} can be reformulated in the following way.
\begin{Lem} 
\label{lllway}
Let $p,q>1$ be coprime integers.
The numbers in the lower left corner of the LLL-diagram are the exponents of 
the terms in $Q_{\{p,q\}}$ with coefficient $1$.
The numbers in the upper right corner are the exponents of the terms in $Q_{\{p,q\}}$ with coefficient $-1$. 
All other coefficients equal $0$.
\end{Lem}

\section{Two proofs of the main (folklore) result}
In terms of inclusion-exclusion polynomials we can reformulate Theorem \ref{basic} as
follows.
\begin{Thm}
\label{basic2}
If $p,q>1$ are coprime integers, then
$P_{S(p,q)}(x)=Q_{\{p,q\}}(x)$. 
\end{Thm}
Our first proof will make use of `what is probably the most versatile tool in numerical semigroup theory' 
\cite[p. 8]{roos}, namely Ap\'ery sets.\\
\noindent {\it First proof of Theorem} \ref{basic2}. The {\it Ap\'ery set} of $S$ with respect to a nonzero $m\in S$ is defined as
$${\rm Ap}(S;m)=\{s\in S:s-m\not\in S\}.$$
Note that
$$S={\rm Ap}(S;m)+m\mathbb Z_{\ge 0}$$
and that ${\rm Ap}(S;m)$ consists of a complete set of residues modulo $m$. Thus we
have
\begin{equation}
\label{hsl}
H_S(x)=\sum_{w\in {\rm Ap}(S;m)}x^w\sum_{i=0}^{\infty}x^{mi}={1\over 1-x^m}\sum_{w\in {\rm Ap}(S;m)}x^w.
\end{equation}
Note that if $S=\langle a_1,\ldots,a_n\rangle$, then ${\rm Ap}(S;a_1)\subseteq \langle a_2,\ldots,a_n\rangle$.
It follows that ${\rm Ap}(S(p,q);p)$ consists of multiples of $q$. The latter set equals the minimal set of multiples of $q$
representing every congruence class modulo $p$ and hence ${\rm Ap}(S(p,q);p)=\{0,q,\ldots,(p-1)q\}$ (see \cite[Proposition 1]{RR} 
or \cite[Example 8.22]{roos}).
Hence
$$H_{S(p,q)}(x)={1+x^q+\cdots+x^{(p-1)q}\over 1-x^p}={1-x^{pq}\over (1-x^q)(1-x^p)}.$$
Using this identity and (\ref{prod4}) easily completes the proof. \qed\\

\noindent {\tt Remark}. This proof is an adaptation of the arguments given in \cite{RR}. Indeed, once we know
the Ap\'ery set of a numerical semigroup $S$, by using \cite[(4)]{RR}, we obtain an expression for
$H_S(x)$ and consequently for $P_S(x)$. Theorem \ref{basic2} is a particular case of \cite[Proposition 2]{RR}, 
with $\{p,q\}=\{a,a+d\}$ and $k=1$.\\

\noindent Our second proof uses the denumerant (see \cite[Chapter 4]{RA} for a survey) 
and the starting point is the observation that
\begin{equation}
\label{flauwflauw}
{1\over (1-x^p)(1-x^q)}=\sum_{j\ge 0}r(j)x^j,
\end{equation}
where $r(j)$ denotes the cardinality of the set
$\{(a,b):a\ge 0,b\ge 0,ap+bq=j\}$. In the terminology of the introduction, we have $r(j)=d(j;p,q)$.
Concerning $r(j)$ we make the following observation.
\begin{Lem}
\label{step1}
Suppose that $k\ge 0$, then $r(k+pq)=r(k)+1$.
\end{Lem}
{\it Proof}. Put $\alpha\equiv kp^{-1}({\rm mod~}q)$, $0\le \alpha<q$
and $\beta\equiv kq^{-1}({\rm mod~}p)$, $0\le \beta<p$ and $k_0=\alpha p+\beta q$.
Note that $k_0<2pq$. We have $k\equiv k_0({\rm mod~}pq)$.
Now if $k\not\in S$, then $k<k_0$ and $k+pq=k_0\in S$ (since $k_0<2pq$).
It follows that if $r(k)=0$, then $r(k+pq)=1$. If $k\in S$, then
$k=k_0+tpq$ for some $t\ge 0$ and we have $r(k)=1+t$, where we use that
$$k=(\alpha+tq)p+\beta q=(\alpha+(t-1)q)p+(\beta+1)p=\cdots=\alpha p +(\beta + tq)p.$$
We see that $r(k+pq)=1+t+1=r(k)+1$. \qed\\

\noindent {\tt Remark}. 
It is not difficult to derive an explicit formula for $r(n)$ (see, e.g., 
\cite[Section 1.3]{Bero} or \cite[pp. 213-214]{niven}). Let $p^{-1},q^{-1}$ denote inverses of $p$ modulo $q$, respectively $q$ modulo $p$.
Then we have
$$r(n)={n\over pq}-\Big\{{p^{-1}n\over q}\Big\}-\Big\{{q^{-1}n\over p}\Big\}+1,$$
where $\{ x \}$ denote the fractional-part function.
Note that Lemma \ref{step1} is a corollary of this formula.\\

\noindent {\it Second proof of Theorem} \ref{basic2}. {}From Lemma \ref{step1} we infer that
\begin{eqnarray}
(1-x^{pq})\sum_{j\ge 0}r(j)x^j & = &\sum_{j=0}^{pq-1}r(j)x^j+\sum_{j=pq}^{\infty}(r(j)-r(j-pq))x^j\nonumber\cr
&=&
\sum_{j=0}^{pq-1}r(j)x^j+\sum_{j\ge pq}x^j=\sum_{j\in S(p,q)}x^j,
\end{eqnarray}
where we used that $r(j)\le 1$ for $j<pq$ and $r(j)\ge 1$ for $j\ge pq$.
Using this identity and (\ref{flauwflauw}) easily completes the proof.\qed

\section{Symmetric numerical semigroups}
A numerical semigroup $S$ is said to be {\it symmetric} if
$$S\cup (F(S)-S) =\mathbb Z,$$
where $F(S)-S=\{F(S)-s|s\in S\}$. Symmetric semigroups occur in the study of monomial curves that are complete intersections,
Gorenstein rings, and the classification of plane algebraic curves, see, e.g. \cite[p. 142]{RA}. 
For example, Herzog and Kunz showed that a Noetherian local ring of dimension one and analytically irreducible is a Gorenstein
ring if and only if its associate value semigroup is symmetric.\\
\indent We will now show
that the selfreciprocity of $Q_{\{p,q\}}(x)$ implies that $S(p,q)$ is symmetric (a well-known result, see, e.g., 
\cite[Corollary 4.7]{roos}).
\begin{Thm}
\label{symmetrie}
Let $S$ be a numerical semigroup. Then $S$ is symmetric if and only if $P_S(x)$ is selfreciprocal.
\end{Thm}
{\it Proof}. If $s\in S\cap (F(S)-S)$, then
$s=F(S)-s_1$ for some $s_1\in S$. This implies that $F(S)\in S$, a contradiction.
Thus $S$ and $F(S)-S$ are disjoint sets. Since every integer $n\ge F(S)+1$ is in $S$ and
every integer $n\le -1$ is in $F(S)-S$, the assertion is equivalent to showing that
\begin{equation}
\label{reversi}
\sum_{0\le j\le F(S)\atop j\in S}x^j+\sum_{0\le j\le F(S)\atop j\in S}x^{F(S)-j}=1+x+\cdots+x^{F(S)},
\end{equation}
if and only if $P_S(x)$ is selfreciprocal. On noting by (\ref{semigppoly}) that
$$x^{F(S)+1}P_S({1\over x})-P_S(x)=1-x^{F(S)+1}+
(x-1)\Big(\sum_{0\le j\le F(S)\atop j\in S}x^j+\sum_{0\le j\le F(S)\atop j\in S}x^{F(S)-j}\Big),$$
we see that $x^{F(S)+1}P_S(1/x)=P_S(x)$ if and only if (\ref{reversi}) holds. Clearly (\ref{reversi}) holds  if and 
only if $S$ is symmetric. \qed\\

\noindent Using the latter result and Theorem \ref{basic2} we infer the following classical fact.
\begin{Thm}
\label{ookeenkeer}
A numerical semigroup of embedding dimension $2$ is symmetric.
\end{Thm}
\noindent Theorem \ref{basic2} together with Theorem \ref{bach} shows that if $e(S)=2$, then
$P_S(x)$ can be written as a product of cyclotomic polynomials. This leads to the following problem.
\begin{Problem}
Characterize the numerical semigroups $S$ for which $P_S(x)$ can
be written as a product of cyclotomic polynomials.
\end{Problem}
Since $P_S(0)\ne 0$, the product cannot involve $\Phi_1(x)=x-1$ and so it is selfreciprocal. Therefore, by
Theorem \ref{symmetrie} such an
$S$ must be symmetric.
Ciolan et al.~\cite{CGM} make some progress towards solving this problem and show, e.g., that $P_S(x)$ can be 
written as a product of cyclotomic polynomials also if $e(S)=3$ and $S$ is symmetric.

\section{Gap distribution}
\label{vijf}
The non-negative integers not in $S$ are called the {\it gaps} of $S$. 
E.g., the gaps in $S(4,7)$ are $1,2,3,5,6,9,10,13$ and $17$.
The number of gaps of $S$ is called the
{\it genus} of $S$, and denoted by $N(S)$. The set of gaps is denoted by $G(S)$. 
The following well-known result holds, cf. \cite[Lemma 7.2.3]{RA} or 
\cite[Corollary 4.7]{roos}.
\begin{Thm}
\label{stern}
We have $2N(S)\ge F(S)+1$ with equality if and only if $S$ is symmetric.
\end{Thm}
{\it Proof}. The proof of Theorem \ref{symmetrie} shows that
$2\#\{0\le j\le F(S): j\in S\}\le F(S)+1$ with equality if and only if $S$ is symmetric.
Now use that $\#\{0\le j\le F(S): j\in S\}=F(S)+1-N(S)$. \qed\\

\noindent {}From (\ref{frob}) and Theorem \ref{basic}
we infer the following well-known result due to Sylvester:
\begin{equation}
\label{sylvester}
F(S(p,q))=pq-p-q.
\end{equation}
{}From Theorem \ref{ookeenkeer}, Theorem \ref{stern} and (\ref{sylvester}),  
we obtain another well-known result of Sylvester:
\begin{equation}
\label{sylvester2}
N(S(p,q))=(p-1)(q-1)/2.
\end{equation}
For four different proofs of (\ref{sylvester}) and more background see \cite[pp. 31-34]{RA}; the shortest proof of (\ref{sylvester}) and
(\ref{sylvester2}) the author knows of is in the book by Wilf \cite[p. 88]{wilf}.\\
\indent Additional information
on the gaps is given by the so-called {\it Sylvester sum}
$$\sigma_k(p,q):=\sum_{s\not\in S(p,q)}s^k.$$
By (\ref{sylvester2}) we have $\sigma_0(p,q)=(p-1)(q-1)/2$.
By (\ref{semigppoly}) and Theorem \ref{basic2} we infer that
\begin{equation}
\label{BB2}
\sum_{j\not\in S(p,q)}x^j={1-Q_{\{p,q\}}(x)\over 1-x}.
\end{equation}
It is not difficult to derive a 
formula for $\sigma_k(p,q)$ for arbitrary $k$. On substituting
$x=e^z$ and recalling the Taylor series expansion
$e^z=\sum_{k\ge 0}z^k/k!$, we obtain from (\ref{BB2})  and (\ref{prod4}) the identity
\begin{equation}
\label{bernie}
\sum_{k=0}^{\infty}\sigma_k(p,q){z^k\over k!}={e^{pqz}-1\over (e^{pz}-1)(e^{qz}-1)}-{1\over e^z-1}.
\end{equation}
We obtain from (\ref{bernie}),  on multiplying by $z$ and using the Taylor series expansion (\ref{generating}), that
$$\sum_{m=1}^{\infty}m\sigma_{m-1}(p,q){z^m\over m!}=\sum_{i=0}^{\infty}B_ip^i{z^i\over i!}\sum_{j=0}^{\infty}B_jq^j{z^j\over j!}
\sum_{k=0}^{\infty}{(pqz)^k\over (k+1)!}-\sum_{m=0}^{\infty}B_m{z^m\over m!}.$$
Equating coefficients of $z^m$ then leads to the following result.
\begin{Thm} {\rm \cite{oystein}.}
\label{stein}
For $m\ge 1$ we have
$$m\sigma_{m-1}(p,q)={1\over m+1}\sum_{i=0}^m\sum_{j=0}^{m-i}{m+1\choose i,j,m+1-i-j} B_iB_jp^{m-j}q^{m-i}-B_m.$$
\end{Thm}
Using this formula we find e.g. that 
$\sigma_1(p,q)={1\over 12}(p-1)(q-1)(2pq-p-q-1)$ (this result is due to Brown and Shiue \cite{brown})
and $\sigma_2(p,q)={1\over 12}(p-1)(q-1)pq(pq-p-q)$. The proof we have given here
of Theorem \ref{stein} is due to R\o dseth \cite{oystein}, with the difference that we gave a different proof
of the identity (\ref{bernie}).\\
\indent By using the formula (\ref{powersum}) for power sums we obtain from Theorem \ref{stein} the identity
$$m\sigma_{m-1}(p,q)=\sum_{r=0}^m {m\choose r}p^{m-r-1}B_{m-r}q^rS_r(p)-B_m,$$
giving rise to the following recursion
formula for $B_m$:
$$B_m={m\over p^m-1}\sigma_{m-1}(p,q)+{q^m\over p(1-p^m)}\sum_{r=0}^{m-1}{m\choose r}\Big({p\over q}\Big)^rB_rS_{m-r}(p).$$
On taking $p=4$ and $q=7$, we obtain the recursion for $B_m$ stated in the introduction.\\
\indent Tuenter \cite{Tuenter} established
the following characterization of the gaps in $S(p,q)$. For every finite function $f$,
$$\sum_{n\not\in S}(f(n+p)-f(n))=\sum_{n=1}^{p-1}(f(nq)-f(n)),$$
where $p$ and $q$ are interchangeable. He shows that by choosing $f$ appropriately one can
recover all earlier results mentioned in this section and in addition the identity
$$\prod_{n\not\in S(p,q)}(n+p)=q^{p-1}\prod_{n\not \in S(p,q)}n.$$
\indent Wang and Wang \cite{WW} obtained results similar to those of Tuenter for the
{\it alternate Sylvester sums} $\sum_{s\not\in S(p,q)}(-1)^ss^k$.

\section{A reproof of Theorem \ref{hoon}}
\label{zes}
As mentioned previously, the gaps for $S(4,7)$ are given by $1,2,3,5,6,9,10,13$ and $17$. One could try to break this down in terms
of {\it gap blocks}, that is blocks of consecutive gaps, (also known in the
literature as {\it deserts} \cite[Definition 16]{FM})): $\{1,2,3\}$, $\{5,6\}$, $\{9,10\}$, $\{13\}$, and $\{17\}$. It
is interesting to compare this with the distribution of the {\it element blocks}, that is finite blocks
of consecutive elements in $S$. 
For $S(4,7)$ we get $\{0\}$, $\{4\}$, $\{7,8\}$, $\{11,12\}$ and $\{14,15,16\}$. The longest
gap block we denote by $g(G(S))$ and the longest element block by $g(S)$.\\
\indent The following result gives some information on gap blocks and element blocks in a numerical semigroup
of embedding dimension $2$.
Recall that the smallest positive integer of $S$ is called the
{\it multiplicity} and denoted by $m(S)$. 
\begin{Lem}
\label{trivialo}${~~}$\\
{\rm 1)} The longest gap block, $g(G(S))$, has length $m(S)-1$.\\
{\rm 2)} The longest element block, $g(S)$, has length not exceeding $m(S)-1$.\\
{\rm 3)} If $S$ is symmetric, then $g(S)=m(S)-1$.
\end{Lem}
{\it Proof}. 1) Let $S=\{s_0,s_1,s_2,s_3,\ldots\}$ be the elements of $S$ written in ascending order, i.e., 
$0=s_0<s_1<s_2<\cdots$. Since $s_0=0$ and $s_1=m(S)$ we have $g(G(S))\ge m(S)-1$. Since all
multiples of $m(S)$ are in $S$, it follows that actually $g(G(S))=m(S)-1$.\\
2) If $g(S)\ge m(S)$, it would imply that we can find $k,k+1,\ldots,k+m(S)-1$ all in $S$ such
that $k+m(S)\not \in S$. This is clearly a contradiction.\\
3) If $S$ is symmetric, then
we clearly have $g(S)=g(G(S))=m(S)-1$.\qed\\
\noindent {\tt Remark}. The second observation was made by my intern Alexandru Ciolan. It allows one to prove Theorem
\ref{ciolan}.\\
\indent Finally, we will generalize a result of Hong et al. \cite{zuidkorea}.
\begin{Thm}
\label{hoon2}
If $p,q>1$ are coprime integers, then  $g(Q_{\{p,q\}}(x))=\min\{p,q\}-1$.
\end{Thm}
{\it Proof}. Note that $g(Q_{\{p,q\}}(x))$ equals the maximum of the longest gap block length and the
longest element block length and hence by Lemma \ref{trivialo} equals $m(S(p,q))-1=\min\{p,q\}-1$.\qed\\
\indent This
result can be easily generalized further.
\begin{Thm}
\label{ciolan}
We have $g(P_S(x))=m(S)-1$.
\end{Thm}
{\it Proof}. Using that $P_S(x)=(1-x)H_S(x)$ and Lemma \ref{trivialo} we infer that 
$g(P_S(x))=\max\{g(S),g(G(S))\}=m(S)-1$. \qed

\section{The LLL-diagram revisited}
It is instructive to indicate (we do this in boldface) the gaps of $S(p,q)$
in the LLL-diagram. They are those elements
$\alpha p + \beta q$ with $0\le \alpha\le q-1$, $0\le \beta\le p-1$
for which $\alpha p + \beta q>pq$. Note that the Frobenius
number equals $(q-1)p+(p-1)q-pq$ and so appears in the top right
hand corner of the LLL-diagram. We will demonstrate this (again) for $p=5$ and $q=7$.
\begin{center}\begin{tabular}{ccc|cccc}
28&33&{\bf 3}&{\bf 8}&{\bf 13}&{\bf 18}&{\bf 23}\\
21&26&31&{\bf 1}&{\bf 6}&{\bf 11}&{\bf 16}\\\hline
14&19&24&29&34&{\bf 4}&{\bf 9}\\
7&12&17&22&27&32&{\bf 2}\\
0&5&10&15&20&25&30
\end{tabular}\end{center}
As a check we can verify that $N(S(p,q))=(p-1)(q-1)/2$ integers appear in boldface. 

On comparing coefficients in the identity $(1-x)\sum_{j\in S(p,q)}x^j=\sum_{j\ge 0}a_{\{p,q\}}(j)x^j$ we
get the following reformulation of Theorem \ref{basic2} at the coefficient level.
\begin{Thm}
\label{hoho2}
If $p,q>1$ are coprime integers, then
$$a_{\{p,q\}}(k)=\cases{1 & if $k\in S(p,q),~k-1\not\in S(p,q)$;\cr
-1 & if $k\not\in S(p,q),~k-1\in S(p,q)$;\cr
0 & otherwise.}$$
\end{Thm}
\begin{cor}
The non-zero coefficients of $Q_{\{p,q\}}$ alternate between $1$ and $-1$.
\end{cor}
\noindent The next result gives an example where an existing result on cyclotomic coefficients yields information
about numerical semigroups.
\begin{Thm}
\label{trivial}
Let $p,q,\rho$ and $\sigma$ be as in Lemma \ref{binary}. If $S=S(p,q)$, then there 
are  $\rho \sigma -1$ gap blocks and $\rho \sigma -1$ element blocks.
\end{Thm}
{\it Proof}. 
In view of Theorem \ref{hoho2} we have $a_{\{p,q\}}(k)=1$ if and only if $k$ is at the start of an element block (including the
infinite block $[F(S)+1,\infty)\cap \mathbb Z$). Moreover, $a_{\{p,q\}}(k)=-1$ if and only if $k$ is at the end of a gap block. 
The proof is now completed by invoking Corollary \ref{carla}. \qed\\

Using Lemma \ref{lllway} and Theorem 
\ref{hoho2} our folklore result can now be reformulated in terms of the LLL-diagram.
\begin{Thm}
\label{hoho3}
Let $p,q>1$ be coprime integers and denote $S(p,q)\cap \{0,\ldots,pq-1\}$ by $T$.
The integers $k\in T$ such that $k-1\not\in T$ are precisely the integers in the 
lower left corner of the LLL-diagram.
The integers $k\not \in T$ such that $k-1\in T$ are precisely the integers in the upper right corner.
If $k$ is not in the lower left or upper right corner, then either $k\in T$ and $k-1\in T$ or
$k\not\in T$ and $k-1\not \in T$.
\end{Thm}
Denote $S(p,q)$ by $S$.
Note that the upper right integer in the lower left corner of the LLL-diagram equals $F(S)+1$ and that
the remaining integers in the lower left corner are all $<F(S)$. This observation together with 
(\ref{sylvester2}) then leads to the following corollary of Theorem \ref{hoho3}.
\begin{cor}
If $p,q>1$ are coprime integers, then
$$\cases{\{0\le k\le F(S):k\in S,k-1\in S\}=(p-1)(q-1)/2-\rho \sigma +1;\cr
\{0\le k\le F(S):k\in S,k-1\not\in S\}=\rho \sigma -1;\cr
\{0\le k\le F(S):k\not\in S,k-1\in S\}=\rho \sigma -1;\cr
\{0\le k\le F(S):k\not\in S,k-1\not\in S\}=(p-1)(q-1)/2-\rho \sigma -1.\cr}$$
\end{cor}
The distribution of the quantity $\rho \sigma$ that appears at various places in this article has been
recently studied using deep results from analytic number theory by Bzd{\k e}ga \cite{Bz} and 
Fouvry \cite{Fouvry}. In particular they are interested in counting the number of
integers $m=pq\le x$ with $p,q$ distinct primes such that $\theta(m)$, the number of non-zero
coefficients of $\Phi_m$, satisfies $\theta(m)\le m^{1/2+\gamma}$, with $\gamma>0$ fixed. (Note
that by Corollary \ref{carla} we have $\theta(m)=2\rho\sigma-1$.)\\

\noindent Acknowledgement. I like to thank Matthias Beck, Scott Chapman, Alexandru Ciolan, Pedro A. Garc\'{\i}a-S\'anchez, Nathan
Kaplan, Bernd Kellner, Jorge Ram\'irez Alfons\'in, Ali Sinan Sertoz, Paul Tegelaar 
and the three referees for helpful comments. Alexandru Ciolan pointed out to me that $g(S)\le m(S)-1$, which 
allows one to prove Theorem \ref{ciolan}.

\medskip\noindent {\footnotesize Max-Planck-Institut f\"ur Mathematik,\\
Vivatsgasse 7, D-53111 Bonn, Germany.\\
e-mail: {\tt moree@mpim-bonn.mpg.de}}
\vskip 5mm

\end{document}